# Proof of the Twin Prime Conjecture
# (Together with the proof of
# Polignac's Conjecture for Cousin Primes)

Marko V. Jankovic

Institute of Electrical Engineering "Nikola Tesla", Belgrade, Serbia

to

勇気

**Abstract** In this paper proof of the twin prime conjecture is going to be presented. Originally very difficult problem (in observational space) has been transformed into a simpler one (in generative space) that can be solved. It will be shown that twin primes could be obtained through two stage sieve process, and that will be used to show that exist infinite number of twin primes. The same approach is used to prove Polignac's conjecture for cousin primes.

## 1 Introduction

In number theory, Polignac's conjecture states: For any positive even number *n*, there are infinitely many prime gaps of size *n*. In other words: there are infinitely many cases of two consecutive prime numbers with the difference *n* [1]. For $n = 2$ it is known as twin prime conjecture.

The results that were presented in the literature and that are the closest to the solution of the twin prime conjecture are the following:

1. Conditioned on the truth of the generalized Elliot-Halberstam Conjecture [2], in [3] it has been shown that there are infinitely many primes' gaps that have value of at most 6.
2. In [3] it has been proved that exists infinitely many primes that have gaps that are not bigger than 246, unconditionally.

In this paper the proof of Polignac's conjecture for $n = 2$ and $n = 4$, is going to be provided. The proof is inspired by the recently proposed proof of Sophie Germain prime conjecture [4]. The problem is addressed in generative space, which means that prime numbers are not going to be analyzed directly, but rather their representatives that are used to produce them. It will be shown that twin primes could be generated by two stage recursion type sieve process. This process will be compared to the other two stage sieve process that leaves infinitely many numbers. The fact that sieve process that generate twin primes leaves more numbers than the other sieve process, will be used to prove that infinitely many twin primes exist.

In the last part of the paper it will be shown that the number of cousin primes is infinite, too.

Without going into the details, here we can say that using a procedure very similar to one proposed in this paper, it is possible to make an elementary proof of Green-Tao theorem [5]. The major difference is that in the case of Green-Tao theorem recursion has the depth that is equal to the length of the arithmetic progression, while the depth of the recursion in the case of twin prime conjecture is 2 (equal to the length of the "arithmetic progression" with two elements).

**Remark 1:** *In this paper any infinite series in the form* $c_1 \cdot l \pm c_2$ *is going to be called a thread defined by number* $c_1$ *(in literature these forms are known as linear factors – however, it seems that the term thread is probably better choice in this context). Here* $c_1$ *and* $c_2$ *are numbers that belong to the set of natural numbers (*$c_2$ *can also be zero and usually is smaller than* $c_1$*) and* $l$ *represents an infinite series of consecutive natural numbers in the form (1, 2, 3, ...).*

## 2 Multiplication tensor

The fundamental theorem of arithmetic states that every integer greater than 1 can be uniquely represented by a product of powers of prime numbers, up to the order of the factors [6]. Having that in mind, an infinite dimensional tensor $\mathbf{M_N}$ that contains all natural numbers only once, is going to be constructed. In order to do that we are going to mark vector that contains all prime numbers with $\mathbf{p}$. So, $p(1) = 2$, $p(2) = 3$, $p(3) = 5$, and so on. Tensor $\mathbf{M_N}$ with elements $m_{i1\ i2\ i3\ \ldots}$ is defined by the following equation ($i_1$, $i_2$, $i_3$, … are natural numbers):

$$m_{i_1 i_2 i_3 \ldots} = p(1)^{i_1-1}\, p(2)^{i_2-1}\, p(3)^{i_3-1} \ldots \quad .$$

The alternative definition is also possible. Now, the following notation is going to be assumed for some infinite size vectors

$$\mathbf{2} = [2^0\ 2^1\ 2^2\ 2^3\ \ldots],\ \mathbf{3} = [3^0\ 3^1\ 3^2\ 3^3\ \ldots],\ \mathbf{5} = [5^0\ 5^1\ 5^2\ 5^3\ \ldots]\ \ldots$$

It is simple to be seen that every vector is marked by bold number that is equal to some prime number and that components of the vector are defined as powers of that prime number, including power zero. Now, the $\mathbf{M_N}$-tensor can be defined as

$$\mathbf{M_N} = \mathbf{2} \circ \mathbf{3} \circ \mathbf{5} \circ \mathbf{7} \circ \ldots,$$

where ○ stands for outer product.

The tensor $\mathbf{M_N}$ is of infinite dimension (equal to number of prime numbers) and size, and contains all natural numbers exactly ones. It is easy to understand why it is so, having in mind the fundamental theorem of arithmetic. We are going to call this type of infinite tensor a half infinite tensor.

Now, we are going to consider presentation of even and odd numbers by infinite dimensional tensor. It is not difficult to be seen that even numbers could be represented by the tensor $\mathbf{M_{NE}}$ whose elements are defined as

$$m_{i_1 i_2 i_3 \ldots} = p(1)^{i_1}\, p(2)^{i_2 - 1}\, p(3)^{i_3 - 1} \ldots \quad ,$$

or as

$$\mathbf{M_{NE}} = \mathbf{2_s} \circ \mathbf{3} \circ \mathbf{5} \circ \mathbf{7} \circ \ldots,$$

where ○ stands for outer product and $\mathbf{2_s}$ is vector defined as

$$\mathbf{2_s} = [2^1\ 2^2\ 2^3\ \ldots].$$

The tensor that represents all odd numbers, $\mathbf{M}_{\mathbf{NO}}$, contains elements defined as

$$m_{i_1 i_2 \ldots} = p(2)^{i_1 - 1} p(3)^{i_2 - 1} \ldots \quad ,$$

or as

$$\mathbf{M}_{\mathbf{NO}} = \mathbf{3} \circ \mathbf{5} \circ \mathbf{7} \circ \ldots,$$

where ○ stands for outer product. This representation is useful, for instance, in the case when we consider implementation of Sundaram or Erathostenes sieves. The tensor $\mathbf{M}_{\mathbf{NO}}$, that represent odd numbers, can also be produced by using shifted tensor $\mathbf{M}_{\mathbf{NE}}$, where elements are defined as

$$m_{i_1 i_2 i_3 \ldots} = -1 + p(1)^{i_1} p(2)^{i_2 - 1} p(3)^{i_3 - 1} \ldots \quad ,$$

or

$$\mathbf{M}_{\mathbf{NO}} = -1 + \mathbf{2}_{\mathbf{s}} \circ \mathbf{3} \circ \mathbf{5} \circ \mathbf{7} \circ \ldots,$$

where ○ stands (as previously defined) for outer product. This tensor is also going to be marked as $\mathbf{M}_{\mathbf{NOS}}$. For instance this representation is good for analysis of the sieve that leaves Mersenne numbers.

## 3 Implementation of sieves

Now, we are going to analyze what is going to happen with $\mathbf{M}_{\mathbf{N}}$-tensor when we implement three different sieves.

First we are going to analyse the Sundaram sieve [7]. Sundaram sieve represents a way to extract prime numbers from natural numbers. Idea is to remove all composite numbers in infinitely many steps using threads.

First, all even numbers (except 2) are removed from the set of natural numbers. Then, it is necessary to remove the composite odd numbers from the rest of the numbers. In order to do that, the formula for the composite odd numbers is going to be analyzed. It is well known that odd numbers bigger than 1, here denoted by $a$, can be represented by the following formula

$$a = 2n + 1,$$

where $n \in N$. It is not difficult to prove that all composite odd numbers $a_c$ can be represented by the following formula

$$a_c = 2(2ij + i + j) + 1 = 2((2j + 1)i + j) + 1. \quad (1)$$

where $i, j \in N$. It is simple to conclude that all composite numbers could be represented by product $(i + 1)(j + 1)$, where $i, j \in N$. If it is checked how that formula looks like for the odd numbers, after simple calculation, equation (1) is obtained. This calculation is presented here. The form $2m + 1$, $m \in N$ will represent odd numbers that are composite. Then the following equation holds

$$2m + 1 = (2i + 1)(2j + 1) \quad ,$$

where $i, j \in N$. Now, it is easy to see that the following equation holds

$$m = 2ij + i + j = (2i + 1)j + i. \quad (2)$$

When all numbers represented by $m$ are removed from the set of odd natural numbers bigger than 1, only the numbers that represent odd prime numbers are going to stay. In other words, only odd numbers that cannot be represented by (1) will stay.

In order to understand the process, a 3D tensor of infinite size, that contains all natural numbers that could be defined by the first 3 primes – 2, 3 and 5 (any three primes could be chosen – it does not change the line of reasoning). That tensor represents 3D sub-tensor of $\mathbf{M}_{\mathbf{N}}$-tensor.

Here we are going to illustrate the frontal, lateral and top slice of this tensor. The frontal, top and lateral slices of the tensor are represented by infinite size 2D matrices given as

| $1$ | $2^1$ | $2^2$ | $2^3$ | ... |
|---|---|---|---|---|
| $3^1$ | $2\cdot 3$ | $2^2\cdot 3$ | $2^3\cdot 3$ | ... |
| $3^2$ | $2\cdot 3^2$ | $2^2\cdot 3^2$ | $2^3\cdot 3^2$ | ... |
| $3^3$ | $2\cdot 3^3$ | $2^2\cdot 3^3$ | $2^3\cdot 3^3$ | ... |
| ... | ... | ... | ... | ... |

Front slice

| $1$ | $2^1$ | $2^2$ | $2^3$ | ... |
|---|---|---|---|---|
| $5^1$ | $2\cdot 5$ | $2^2\cdot 5$ | $2^3\cdot 5$ | ... |
| $5^2$ | $2\cdot 5^2$ | $2^2\cdot 5^2$ | $2^3\cdot 5^2$ | ... |
| $5^3$ | $2\cdot 5^3$ | $2^2\cdot 5^3$ | $2^3\cdot 5^3$ | ... |
| ... | ... | ... | ... | ... |

Top slice

| $1$ | $5^1$ | $5^2$ | $5^3$ | ... |
|---|---|---|---|---|
| $3^1$ | $5\cdot 3$ | $5^2\cdot 3$ | $5^3\cdot 3$ | ... |
| $3^2$ | $5\cdot 3^2$ | $5^2\cdot 3^2$ | $5^3\cdot 3^2$ | ... |
| $3^3$ | $5\cdot 3^3$ | $5^2\cdot 3^3$ | $5^3\cdot 3^3$ | ... |
| ... | ... | ... | ... | ... |

Lateral slice

Implementation of Sundaram sieve is presented in Figure 1.

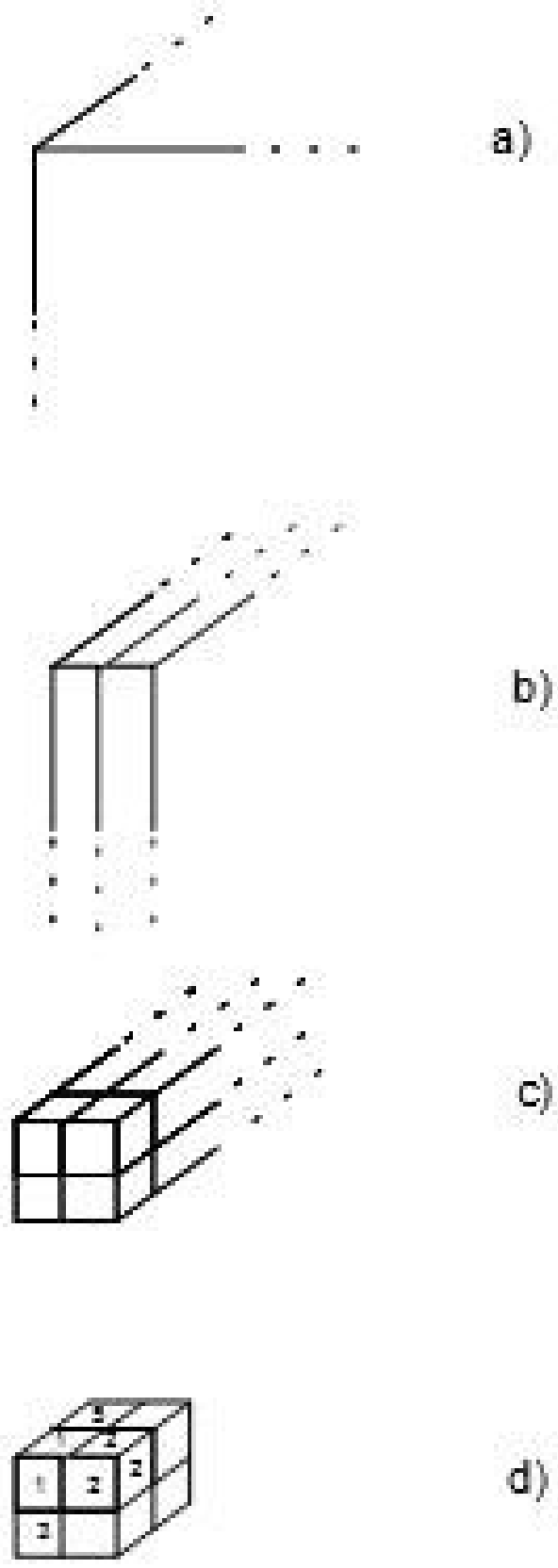

Figure 1. a) 3D infinite tensor; b) collapse of the tensor after removal of the thread defined by p(1); c) additional collapse of the tensor after removal of the thread defined by p(2); final collapse of the tensor after removal of the thread defined by p(3).

From Fig1. it can be seen that removal of a single thread defined by first prime number will result in collapse of the tensor along the dimension that is defined by that prime number. The whole sub-tensor is going to collapse to one element that contains prime number that defines the corresponding thread and that dimension of the tensor. All other elements along that dimension are empty (the empty elements are introduced in order to have reasonable tensor representation). At the end of the process we are going to have 2x2x2 tensor that have only 4 elements that are non-empty. One is top, left corner that is filled by 1, and 3 of his neighbors along the orthogonal dimensions that are filled with first three prime numbers. All other elements (4 in total) are empty. In general case only those elements of the tensor that are going to be nonempty after implementation of full Sundaram sieve, are given by the following equation (tensor will have 2x2x2x … dimension)

$$m_{111...}=1\,;m_{12111...}=p(1);m_{11211...}=p(2);m_{1112 1...}=p(3);...m_{111...2(\text{on position}\,k+1)11...}=p(k);... \quad .$$

The second sieve that is going to be considered is a sieve that is obtained by simple extension of Sundaram sieve – when you remove one thread defined by some prime number, you also remove a prime number that defines that thread. We will call this sieve extended Sundaram sieve or eSuS sieve.

The third sieve of interest is sieve that is the one which we are going to call Mersenne sieve. In this sieve, in first step all even numbers are removed. Then, from the odd numbers are removed the threads defined by the following equation:

$$a_i=2\,p(i)\,j-1,\;\; i=2,3,4,... \qquad (3)$$

where $j \in N$. It is not difficult to be seen that only numbers that are going to be left are Mersenne numbers, or numbers in the form

$$a=2^i-1, i=1,2,3,...$$

This can be easily concluded from the $\mathbf{M}_{\mathbf{NOS}}$ representation of odd numbers – removal of individual threads causes the collapse of the $\mathbf{M}_{\mathbf{NOS}}$ tensor along that dimension, and only dimension that is not going to be affected is the one that is defined by prime number 2.

Now, we are going to compare three sieves that we mentioned in this section. We are going to compare their characteristics and the number of numbers left (that are smaller than some natural number $n$) after implementation of several steps of the sieve.

Before we start with comparison, we are going to analyse two elementary experiments with enumerated balls and boxes (sieves can be presented in a similar context). As it is going to be shown in the following examples, if you stick to the numbers on the balls that are moved, rather than number of the balls that are moved, numbers can, even, create a small problem.

First experiment: Imagine that we have infinite number of balls with all natural numbers written on them exactly once that are placed in the source box (SB), that has the size equal to the number of natural numbers, and that we have another, experimental box (EB) of proper size.

In the moment 1 minute before midnight, we are going to move balls with numbers 1 to 10 on them from SB to EB, and remove from the EB the ball with number 10 on it. In the moment ½ minute before midnight, balls with numbers 11-20 are transferred from SB to EB, and ball with number 20 on it is removed from EB. We continue the process at the moments $1/2^n$, $n \in N$ minute before midnight – transfer the balls with numbers form $n*10+1$ to $(n+1)*10$ from SB to EB, and remove ball with number $(n+1)*10$ on it from EB.

Now, we can try to answer the following question: What is the number of the balls in EB at midnight? The answer is obvious and everybody will answer that that number is infinite.

Second experiment: Imagine, again, that we have infinite number of balls with all natural numbers written on them exactly once that are placed in the source box (SB), that has the size equal to the number of natural numbers, and that we have another, experimental box (EB) of proper size.

In the moment 1 minute before midnight, we are going to move balls with numbers 1 to 10 on them from SB to EB, and remove from the EB the ball with number 1 on it. In the moment ½ minute before midnight, balls with numbers 11-20 are transferred from SB to EB, and the ball with number 2 on it is removed from EB. We continue the process at the moments $1/2^n$, $n \in N$ minute before midnight – transfer the balls with numbers from $n$*10+1 to ($n$+1)*10 from SB to EB, and remove ball with number ($n$+1) on it from EB.

Now, we can try, again, to answer the following question: What is the number of the balls in EB at midnight? Again, the answer is obvious and everybody will answer that that number is infinite. However, if you are asked to give an example of the ball with any specific number on it, that is still in the EB, you will not be able to do it. The reason is quite obvious – for any number you choose, you can specify a moment in time in which the ball with that number on it has been removed form the EB. This process represents an algorithm how you can remove first 1/10 of the natural numbers from the set of natural numbers. What can be interesting to notice is that we can actually define shallow and deep infinity and shallow and deep limes – in this case shallow limes of the number of the balls in the EB is 0, while deep limes is infinity. This can be useful in analysing some sieves where potentially it is easy to prove some properties in shallow infinity while, at the same time, it is difficult to be done in deep infinity. However, this discussion is out of the scope of this paper.

So, numbers on the balls are not relevant for reasoning – in both experiments that were previously analysed, the process was following – put 10 balls in the EB and then remove one – or very simplified, put nine balls in the EB in every step. If there is no collapse of elementary reasoning (CER), we can safely conclude that number of the balls in the EB at midnight is infinite. Actually, it is not difficult to be seen that two previously mentioned experiments are the special cases of a more general experiment in which in every step ten balls are put in the EB and one of the existing balls in the EB is removed completely randomly.

Now, comparison of Sundaram (Erathostenes), eSuS and Mersenne sieves is going to be preformed. In the text that follows Sundaram sieve will be marked as $SuS^0$ and Mersenne sieve as MeS. The sieves are going to be compared in three categories:

1. how many numbers smaller than some natural number $n$ are they going to leave,

2. how many threads are they going to need to achieve this,

3. if some of the numbers left, are going to be removed by the threads that will enable to check number of numbers smaller than some natural number $m > n$ that are left by implementation of sieve in some future search.

The following table contains comparison of the three sieves of interest. We denote with $\pi(n)$ number of primes smaller than $n$.

Table 1. Comparison of Sundaram (Erathostenes), extended Sundaram and Meresenne sieves

| | $SuS^0$ | eSuS | MeS |
|---|---|---|---|
| Number of numbers smaller than $n$ left | $\pi(n)+1\approx\pi(n)$ | $\pi(n)-\pi(\sqrt{n})+1\approx\pi(n)$ | $\text{floor}(\log_2(n))$ |
| Number of threads used | $\pi(\sqrt{n})$ | $\pi(\sqrt{n})$ | $\pi(n/2)$ |
| Will some numbers are going to be removed in next steps | no | yes | no |

As it is indicated in the table, in the text that follows we are going to ignore that sieves $SuS^0$ and eSuS will leave number 1. Also, we are going to consider that the value $\pi(\sqrt{n})$ can be ignored in comparison to value $\pi(n)$ for $n$ that is large enough (and that is for all cases of interest).

Here, the following inequalities are going to be presented without presenting the proofs (although the proofs are not difficult):

$$n-\sqrt{n}\gg\pi(\sqrt{n}) \quad \text{, for some } n \text{ big enough} \tag{4}$$

$$\pi\left(\frac{\pi(n)}{2}\right)>\pi(\sqrt{n}) \quad \text{, for some } n \text{ that is big enough.} \tag{5}$$

The value of $n$ in (5), that is big enough, can be precisely found by using the fact that prime counting function is non-decreasing and that square root function is strictly increasing, but here that is not of particular interest (see Fig 2.). In equation (4) what is the $n$ that is big enough is something that is defined by the context in which inequality is used.

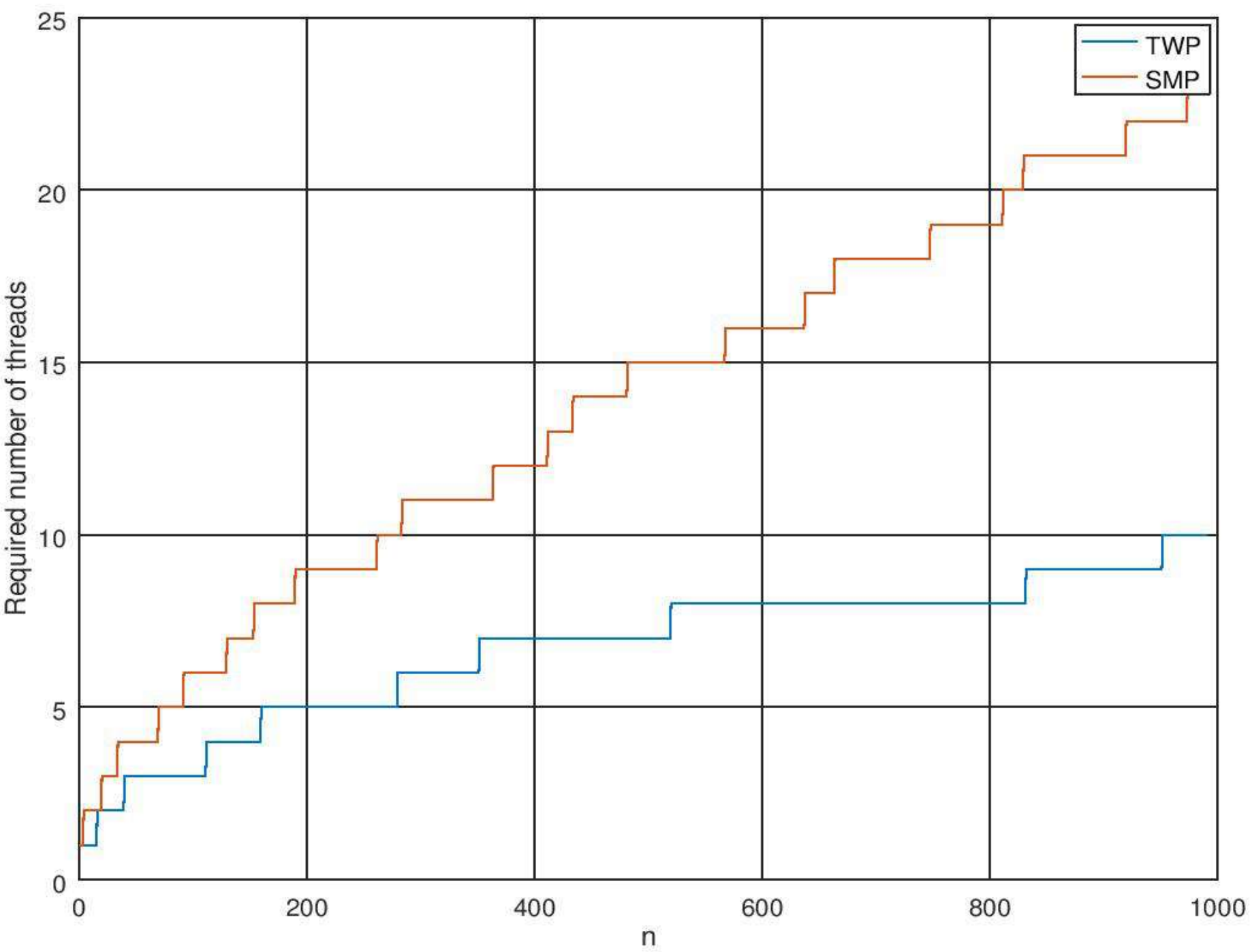

Fig. 2. Comparison of the number of threads required for the realization of the MeS and $SuS^{-1}$ sieves for generation of semi Mersenne primes (SMP) and twin primes (TWP) numbers smaller than some natural number *n*

Here, an additional difference between $SuS^0$ and MeS on one side and eSuS on the other side is going to be stressed. It is not difficult to be understood that after infinite number of steps, or removal of infinite number of threads, we are not going to be able to specify a single number left by eSuS although the number of numbers left is infinite. This is another example in which limes in shallow infinity is zero and limes in deep infinity is infinite. For $SuS^0$ and MeS both shallow and deep infinity limeses are infinite. Also one more thing is worth noticing: $SuS^0$ will leave infinitely more numbers left than the eSuS (every thread leaves one number more), however that infinite difference is very small comparing to the overall number of numbers left (asymptotically the ratio of those two numbers is zero).

## 4 Proof of the twin prime conjecture

It is well known that every two consecutive odd numbers ($ps_k$, $pl_k$) between two consecutive odd numbers divisible by 3 (e.g. 9 11 13 15, or 39 41 43 45), can be expressed as

$$ps_k = 6k - 1,\ \ pl_k = 6k + 1,\ k \in N. \tag{6}$$

Twin prime numbers are obtained in the case when both $ps_k$ and $pl_k$ are prime numbers. If any of the $ps_k$ or $pl_k$ (or both) is a composite number, then we cannot have twin primes. In the text that follows we will call numbers $ps_k$ - numbers in *mps* form and numbers $pl_k$ - numbers in *mpl* form.

Here we are going to present a two stage sieve process that can be used for generation of twin primes. In the first stage we are going to produce prime numbers by removing all composite numbers from the set of natural numbers. In the second stage, we are going to remove all prime numbers that have a bigger odd neighbor that is a composite number. At the end, only the prime numbers in the *mps* form, that represent the smaller number of a twin prime pair, are going to stay. Their number is equal to the half of the number of twin primes. It is going to be shown that that number is infinite. It is easy to check that all numbers in *mpl* form are going to be removed from the set, since their bigger odd neighbors are composite numbers divisible by 3. (Of course it would be possible to organize stage 2 in such a way to remove all odd numbers that have smaller neighbor that is composite. In that case what would be left are prime numbers in *mpl* form, that represent bigger primes in twin pairs, and all primes in *mps* form would be removed, since they have a smaller neighbor divisible by 3. This is completely equivalent process to the first one and will not be analyzed further.)

**STAGE 1**

Prime numbers can be obtained by implementation of Sundaram sieve as it is explained in the previous section.

When all numbers represented by $m$ in (2) are removed from the set of odd natural numbers bigger than 1, only the numbers that represent odd prime numbers are going to stay. In other words, only odd numbers that cannot be represented by (1) will stay. As it was already said, sieve defined by (2) is going to be marked as $SuS^0$.

The numbers that are left after this stage are prime numbers, and their number is $\pi(n)$. From [8] we know that the following holds

$$\pi(n) \approx \frac{n}{\ln(n)} \quad .$$

**STAGE 2**

What was left after the first stage are prime numbers. With the exception of number 2, all other prime numbers are odd numbers. All odd primes can be expressed in the form $2n + 1$, $n \in N$. It is simple to understand that their bigger odd neighbor must be in the form $2n + 3$, $n \in N$. Now, we should implement a second step in which we are going to remove number 2 (since 2 cannot make twin pair) and all odd primes (in the form $2m+1$) whose bigger neighbor is composite number in the form $2m + 3$, $m \in N$. If we make the same analysis again, it is simple to understand that $m$ must be in the form

$$m = 2ij + i + j - 1 = (2i + 1)\,j + i - 1. \tag{7}$$

All numbers (in observational space) that are going to stay must be numbers in *mps* form and they represent a smaller primes of the twin pairs (it is simple to understand that prime numbers in *mpl* form have neighbors that are composite odd numbers divisible by 3). So, their number is equal to the half of the number of twin primes. The sieve defined by (7) is going to marked as $SuS^{-1}$.

It is not difficult to understand (prove) that sieves $SuS^0$ and $SuS^{-1}$ can be implemented in reverse order and the final result will be the same. Those two sieves are almost equivalent, if the number of numbers left after the implementation of the sieve are considered. In the case when $SuS^{-1}$ is implemented first, the number of numbers smaller than natural number $n$, that are left after implementation of sieve is going to be marked as $\pi_{-1}(n)$. It is very simple to prove that the following equation holds

$$\pi(n) - \pi_{-1}(n) \leq 1.$$

Let us mark the number of twin primes with $\pi_{G2}$. Also, let us define the number of numbers that is left after consecutive implementations of Sundaram sieve, $SuS^0$, and Mersenne sieve, MeS, as $p_{SMP}$. The numbers obtained after implementation of those sieves are going to be called semi Merssenne primes, or SMP. The second stage MeS sieve is applied on prime number indexes. In that case it is easy to understand that the following equation would hold ($n \in N$)

$$p_{SMP}(n) = floor(\log_2(\pi(n))) \quad , \tag{8}$$

where $p_{SMP}(n)$ represents the number of SMP smaller than some natural number $n$. Since the twin primes are obtained by implementation of the $SuS^0$ in the first stage and sieve that is similar to it, $SuS^{-1}$ in the second stage, it is not difficult to conclude that numbers $\pi_{G2}$ is bigger than the $p_{SMP}$. In order to understand why it is so, we are going to analyze (3) and (7) in more detail.

It is not difficult to be seen that $m$ in (3) and (7) is represented by the threads that are defined by odd prime numbers. For details see Appendix A. Now we are going to compare sieves $SuS^{-1}$ and MeS step by step, for a few initial steps (analysis can be easily extended to any number of steps). Removal of number 2 in the second stage of $SuS^{-1}$ is ignored. What we would like to say is that the values of the fractions presented in the table are asymptotically correct, but in the finite case they are only approximately correct.

**Table 2 Comparison of the stages 2 for generation of semi-Mersenne primes and generation of twin primes – threads defined by a few smallest primes**

| Step | Stage 2 - MeS | Step | Stage 2 - $SuS^{-1}$ |
|---|---|---|---|
| 1 | Remove even numbers<br><br>amount of numbers left is 1/2 | 1 | Remove numbers defined by thread defined by 3 (obtained for $i$ = 1) amount of numbers left is 1/2 |
| 2 | Remove numbers defined by thread defined by 3 (obtained for $i$ = 1) amount of numbers left is 2/3 of the numbers that are left after previous step | 2 | Remove numbers defined by thread defined by 5 (obtained for $i$ = 2) amount of numbers left is 3/4 of the numbers that are left after previous step |
| 3 | Remove numbers defined by thread defined by 5 (obtained for $i$ = 2) amount of numbers left is 4/5 of the numbers that are left after previous step | 3 | Remove numbers defined by thread defined by 7 (obtained for $i$ = 3) amount of numbers left is 5/6 of the numbers that are left after previous step |
| 4 | Remove numbers defined by thread defined by 7 (obtained for $i$ = 3) amount of numbers left is 6/7 of the numbers that are left after previous step | 4 | Remove numbers defined by thread defined by 11 (obtained for $i$ = 5) amount of numbers left is 9/10 of the numbers that are left after previous step |
| 5 | Remove numbers defined by thread defined by 11 (obtained for $i$ = 5) amount of numbers left is 10/11 of the numbers that are left after previous step | 5 | Remove numbers defined by thread defined by 13 (obtained for $i$ = 6) amount of numbers left is 11/12 of the numbers that are left after previous step |

From the table, it can be noticed that threads defined by the same prime number in the first and the second column will not remove the same percentage of numbers. The reason is obvious – consider for instance the thread defined by 3: in the first column it will remove 1/3 of the numbers left, but in the second column it will remove ½ of the numbers left, since the thread defined by 3 in stage 1 has already removed one third of the numbers (odd numbers divisible by 3 in observation space). So, only odd numbers (in observational space) that give residual 1 and -1 when they are divided by 3 are left, and there is approximately same number of prime numbers that give residual -1 and numbers that give residual 1, when the number is divided by 3. Same way of reasoning can be applied for all other threads defined by the same prime in different columns that represent different second stage sieves - MeS and $SuS^{-1}$. More rigorous poof for the values of the fractions in the second column (and all others that are not presented in the table), based on Dirichlet's theorem on arithmetic progressions [9], is presented in Appendix B.

From Table 2 can be seen that in every step, except step 1, threads in the second column will leave bigger percentage of numbers than the corresponding threads in the first column. This could be easily understood form the analysis that follows:

- suppose that we have two natural numbers $j$, $k$ such that $j - 1 \geq k$ ($j, k \in N$), then the following set of equations is trivially true

$$j+k-1 \geqslant 2k$$

$$-j-k+1 \leqslant -2k$$

$$jk-j-k+1 \leqslant jk-2k$$

$$(j-1)(k-1) \leqslant (j-2)k$$

$$\frac{k-1}{k} \leqslant \frac{j-2}{j-1}$$

The equality sign holds only in the case $j = k + 1$. In the set of prime numbers there is only one case when $j = k + 1$ and that is in the case of primes of 2 and 3. In all other cases $p(i) - p(i-1) > 1$, ($i > 1$, $i \in N$, $p(i)$ is $i$-th prime number). So, in all cases $i > 2$

$$\frac{p(i-1)-1}{p(i-1)} < \frac{p(i)-2}{p(i)-1} \quad .$$

From Table 2 (or last equation) we can see that bigger number of numbers is left in every step of stage 2 then in the stage 1 (except 1st step). From that, we can conclude that after every step bigger than 1, part of the numbers that is left in stage 2 after implementation of $SuS^{-1}$ is bigger than the number of numbers left in the stage 2 after implementation of MeS (that is also true if we consider amount of numbers left after removal of all numbers generated by threads that are defined by all prime numbers smaller than some natural number).

Having in mind previous analysis, it can be safely concluded that the following equation holds

$$\frac{\pi_{G2}}{2} > p_{SMP} = \lim_{n \to \infty} p_{SMP}(n) \quad .$$

From previous inequality it can be concluded that the following equation must hold

$$\pi_{G2} > p_{SMP} = \lim_{n \to \infty} p_{SMP}(n) \quad .$$

Having in mind (8), and since it it easy to show that following holds

$$\lim_{n \to \infty} \log_2(\pi(n)) = \infty \quad ,$$

then it is easy to understand that the following equation holds

$$p_{SMP} = \lim_{n \to \infty} p_{SMP}(n) = \infty \, .$$

Now, we can safely conclude that the number of twin primes is infinite. That concludes the proof.

## 5. Estimation of the number of twin primes

Here we will state the following conjecture: for $n$ big enough, number of twin primes smaller than some natural number $n$, marked as $\pi_{G2}(n)$ is given by the following equation

$$\pi_{G2}(n) \sim 4\, C_2 \pi(\pi(n)) \approx 4\, C_2 \frac{\pi(n)}{\ln(\pi(n))} ,$$

where $C_2$ represents the twin prime constant [10]. Why it is reasonable to make such conjecture is explained in Appendix C. More precisely, instead of value $\pi\,(\pi(n))$ in the previous formula should stay value $\pi(\pi(n)) - \pi(\sqrt{\pi(n)}) + 1$. However, based on the discussion after the Table 1, for $n$ big enough value $-\pi\,(\sqrt{\pi(n)}) + 1$ can be ignored.

If we mark the number of primes smaller than some natural number $n$ with $\pi(n) \approx f\,(n)$, where function $f\,(n)$ gives good estimation of the number of primes smaller than $n$, than $\pi_{G2}(n)$, for $n$ big enough, is given by the following equation

$$\pi_{G2}(n) \sim 4C_2 \cdot f(f(n)).$$

In particular case $f(n) = Li(n)$, the following equation holds

$$\pi_{G2}(n) \sim 4C_2 \cdot \int_2^n \left( \frac{dx}{\ln(x)\ln\left(\int_2^x \left(\frac{dt}{\ln(t)}\right)\right)} \right).$$

For numbers $n$ in the range $[10^3 , 10^{16}]$, the following equation gives good estimation of the number of twin primes smaller than $n$

$$\pi_{G2}(n) \approx 0.985 \cdot \left( \frac{\frac{2}{3}\log(n)}{\frac{2}{3}\log(n)+1} \right) \cdot 4C_2 \cdot \int_2^n \left( \frac{dx}{\ln(x)\ln\left(\int_2^x \left(\frac{dt}{\ln(t)}\right)\right)} \right).$$

In particular case of recently proposed $f(n) = MoLi(n)$ [11], the following equation holds

$$\pi_{G2}(n) \sim 4C_2 \cdot \int_2^n \left( \left( \frac{1}{\ln(x+\sqrt{x})} - \frac{1}{2} \int_2^x \left( \frac{dt}{(t\sqrt{x}+x)(\ln^2(t+\sqrt{x}))} \right) \right) \frac{dx}{\ln\left( \int_2^x \left( \frac{dt}{\ln(t+\sqrt{x})} \right) + \sqrt{n} \right)} \right).$$

For numbers $n$ in the range $[10^2 , 10^6]$, the following equation gives good estimation of the number of twin primes smaller than $n$

$$\pi_{G2}(n) \sim 2 \int_2^n \left( \left( \frac{1}{\ln(x+\sqrt{x})} - \frac{1}{2} \int_2^x \left( \frac{dt}{(t\sqrt{x}+x)(\ln^2(t+\sqrt{x}))} \right) \right) \frac{dx}{\ln\left( \int_2^x \left( \frac{dt}{\ln(t+\sqrt{x})} \right) + \sqrt{n} \right)} \right).$$

Another formula for the estimation of the number of twin primes smaller than some natural number $n$, can be derived from the proof presented in this paper. If the function $f(n) = MoLi(n)$ is considered as good approximation function for estimation of the number of primes smaller than some natural number $n$, the following equation holds for the estimation of the number of twin primes smaller than some natural number $n$

$$\pi_{G2}(n) \sim 4C_2 \cdot \int_{2}^{\pi(n)} \frac{dx}{\ln\left(x+\sqrt{\pi(n)}\right)} \quad .$$

For small values of $n$, it is necessary to use correction term (that is smaller than 1, but asymptotically tends to 1 as $n$ tends toward infinity). Selection of this term will not be analyzed here.

## 6. Proof that the number of cousin primes is infinite

The cousin primes are successive prime numbers with gap 4. It is clear that cousin primes represent pairs of odd numbers that surround odd number divisible by 3 (e.g. (7 9 11), or (13 15 17)). A pair can only represent a cousin primes if both those numbers are primes. So, if we denote a pair of odd numbers that surround an odd number divisible by 3 as $pl_k = 6k + 1$ and $ps_k = 6(k + 1) - 1$, $k \in N$, these numbers can represent cousin primes only in the case when both $pl_k$ and $ps_k$ are prime numbers. If any of the $ps_k$ or $pl_k$ (or both) is a composite number, then we cannot have cousin primes.

Here, similar to the case of twin primes we are going to create a two stage process for generation of cousin primes.

**STAGE 1**

Using the same methodology as previously, generate all prime numbers. In order to do that, from the set of all natural numbers bigger than 1, remove all even numbers (except 2) and all odd numbers generated by equation (2).

**STAGE 2**

What was left after first stage are prime numbers. With the exception of number 2, all other prime numbers are odd numbers. All odd primes can be expressed in the form $2n + 1$, $n \in N$. It is simple to understand that their bigger odd cousin must be in the form $2n + 5$, $n \in N$. Now, we should implement a second step in which we are going to remove number 2 (since 2 cannot make cousin pair) and all odd primes (in the form $2m + 1$) that have bigger odd composite cousin in the form $2m + 5$, $m \in N$. If we make the same analysis like in the case of twin primes, it is simple to understand that $m$ must be in the form

$$m = 2ij + i + j - 2 = (2i + 1)j + i - 2. \tag{9}$$

All numbers (in observational space)| that are going to stay must be numbers in *mpl* form and they represent a smaller primes of the cousin pairs (it is simple to understand that prime numbers in *mps* form have cousins that are composite odd numbers divisible by 3). Their number is equal to the half of cousin primes.

Now, using the same method like in the case of the twin prime conjecture, it can be proved that exists infinitely many cousin primes.

Let us mark the number of cousin primes smaller than some natural number $n$ with $\pi_{G4}(n)$. Here we will state the following conjecture (see Appendix C): for $n$ big enough, number of cousin primes is given by the following equation

$$\pi_{G4}(n) \sim 4C_2 \frac{\pi(n)}{\ln(\pi(n))}.$$

If we mark the number of primes smaller than some natural number $n$ with $\pi(n) \approx f(n)$, where function $f(n)$ gives good estimation of the number of primes smaller than $n$, than $\pi_{G4}(n)$, for $n$ big enough, is given by the following equation

$$\pi_{G4}(n) \sim 4C_2 \cdot f(f(n)).$$

*Note: Here we can see that constant $C_2$ has a misleading name. It is connected with repeated (recursive) implementation of a sieve that produces certain type of prime numbers (in the second step same sieve is applied on the set depleted by the first sieve) which is also, but not exclusively, connected to the twin primes. It seems that better notation for that constant would be ${}^2S$.*

**APPENDIX A.**

Here it is going to be shown that $m$ in (2) is represented by threads defined by odd prime numbers. Now, the form of (2) for some values of $i$ will be checked.

***<u>Case i = 1:</u>*** $m = 3j + 1$,

***<u>Case i = 2:</u>*** $m = 5j + 2$,

***<u>Case i = 3:</u>*** $m = 7j + 3$,

***<u>Case i = 4:</u>*** $m = 9j + 4 = 3(3j + 1) + 1$,

***<u>Case i = 5:</u>*** $m = 11j + 5$,

***<u>Case i = 6:</u>*** $m = 13j + 6$ ,

***<u>Case i = 7:</u>*** $m = 15y + 7 = 5(3j + 1) + 2$,

***<u>Case i = 8:</u>*** $m = 17j + 8$,

It can be seen that $m$ is represented by the threads that are defined by odd prime numbers. From examples (cases $i = 4$, $i = 7$), it can be seen that if $(2i + 1)$ represent a composite number, $m$ that is represented by thread defined by that number also has a representation by the thread defined by one of the prime factors of that composite number. That can be proved easily in the general case, by direct calculation, using representations similar to (2). Here, that is going to be analyzed. Assumimg that $2i + 1$ is a composite number, the following holds

$$2i + 1 = (2l + 1)(2s + 1)$$

where ($l, s \in N$). That leads to

$$i = 2ls + l + s.$$

The simple calculation leads to

$$m = (2l + 1)(2s + 1)j + 2ls + l + s = (2l + 1)(2s+1)j + s(2l + 1) + l$$

or

$$m = (2l+1)((2s+1)j + s) + l$$

which means

$$m = (2l + 1)f + l,$$

and that represents the already exiting form of the representation of $m$ for the factor $(2l + 1)$, where

$$f = (2s + 1)j + s.$$

In the same way this can be proved for (3), (7) and (9) .

Note: *It is not difficult to understand that after implementation of stage 1, the number of numbers in residual classes of some specific prime number are equal. In other words, after implementation of stage 1, for example, all numbers divisible by 3 (except 3, but it does not affect the analysis) are removed. However, the number of numbers in the forms 3k + 1 and 3k + 2 (alternatively, 3k – 1) are equal. The reason is that the thread defined by any other prime number (bigger than 2) will remove the same number of numbers from the numbers in the form 3k + 1 and from the numbers in the form 3k + 2. It is simple to understand that, for instance, thread defined by number 5, is going to remove 1/5 of the numbers in form 3k + 1 and 1/5 of the numbers in form 3k + 2. This can be proved by elementary calculation. That will hold for all other primes and for all other residual classes.*

**APPENDIX B.**

Now we are going to show that the inputs in the second column of Table 2 are correct. In order to do that, we are going to use Dirichlet's theorem on arithmetic progressions [9]. The theorem states that for any two positive coprime integers $a$ and $d$, there are infinitely many prime numbers in the form $a + nd$, where $n$ is also positive integer. Beside that, theorem also proves that for a given value of $d$, proportion of primes in each of progressions $a + nd$, asymptotically, is $1/\varphi(d)$, where $\varphi(d)$ represents Euler's totient function [12] that represents number of feasible progression for a given $d$, such that $a$ and $d$ are coprimes.

In the analysis that follows $k$ represent natural number, while $n$ represents nonnegative integer numbers. Also, in order to simplify analysis, it is assumed that reader is capable to understand when certain context requires use of only odd or only even numbers $n$ and/or $k$.

It is easy to understand that any thread in generative space, defined by some prime number in (7), will generate the thread in observational space that is defined by the same prime number, but with a different residual class. For instance, thread $3k$ in generative space will produce the thread $3(2k) +1$ in observational space and so on. So, from now on, we are going to analyze numbers in observational space in order to make analysis easier.

Now, in the Step 1 of Stage 2, for generation of twin primes, the numbers that are going to be removed are generated by the thread defined by the prime number 3. That corresponds to the thread $3(2k)+1$ in observational space. In that case one half of the prime numbers are going to be removed. That follows directly from the Dirichlet's theorem [9], since all prime numbers can be expressed only in the form $3n+1$ or $3n+2$.

In the next step we are going to analyze what is going to happen when we remove thread defined by number 5, and which is given by $5(2k) + 3$ in observational space. In order to understand that, we are going to represent all numbers by the 15 threads defined by number 15. Those threads are defined by the following progressions:

$15n$+1, $15n$+2, $15n$+3, $15n$+4, $15n$+5, $15n$+6 $15n$+7, $15n$+8, $15n$+9, $15n$+10, $15n$+11, $15n$+12, $15n$+13, $15n$+14 and $15n$+15.

We know that in the first stage (generation of prime numbers) numbers divisible by 3 (1/3 of the threads) defined by threads $15n$+3, $15n$+6, $15n$+9, $15n$+12 and $15n$+15 are going to be removed, as well as numbers divisible by 5 (1/5 of the threads left) defined by the threads $15n$+5 and $15n$+10.

The threads that are left are:

$15n$+1, $15n$+2, $15n$+4, $15n$+7, $15n$+8, $15n$+11, $15n$+13 and $15n$+14.

Based on Dirichlet's theorem we know that each of these threads contain 1/8 of the prime numbers.

In the first step of second stage odd numbers in the form $3k + 1$, should be removed. That means that half of the threads that can be represented in the form $3k$ +1 are going to be removed. Those threads are $15n$+1, $15n$+4, $15n$+7 ad $15n$+13.

So, the threads that are left are $15n$+2, $15n$+8, $15n$+11 and $15n$+14, and each of them contains 1/8 of the prime numbers.

In the second step of the second stage, numbers defined by thread defined by 5, in the form $5(2k)$+3 have to be removed. The only thread that is left and that can be expressed in the form $5(2k)$+3 is thread $15n$+8. Now, we can easily conclude that number of primes that is removed by the thread defined by number 5 is ¼ of the numbers left (or that ¾ of the available numbers is left).

However, it is difficult to generalize the proposed method for the other steps in Stage 2. So, an alternative method is going to be analyzed.

After all numbers in the form $3(2k) + 1$ are removed, we know that all odd prime numbers that are left have to be in the form $3(2k+1) + 2$, or, for the sake of simplicity, in the form $3k + 2$ (and reader should have in mind that we are talking only about odd numbers, since all even numbers were removed in the first step of Stage 1). We know that all numbers that are left have to be in some of

the following forms

$3(5n+1)+2$, $3(5n+2)+2$, $3(5n+3)+2$, $3(5n+4)+2$ and $3(5n+5)+2$, or

$$15n + 3*1 + 2\ ,\ 15n + 3*2 + 2,\ 15n + 3*3 + 2\ ,\ 15n + 3*4 + 2\ ,\ 15n + 3*5 + 3*0 + 2. \qquad (B.1)$$

Since all forms in (B.1) contain the term(s) divisible by 15 (and consequently divisible by 5), it is clear that additional forms that are going to be removed, will be removed based on the analysis of the following expressions

$$3*1 + 2\ ,\ 3*2 + 2,\ 3*3 + 2\ ,\ 3*4 + 2\ ,\ 3*0 + 2\ . \qquad (B.2)$$

We know that in the first stage thread that is divisible by 5 has to be removed and in the second step of the Stage 2, thread that is in the form $5k + 3$, has to be removed. We can see that all five terms in equation (B.2) represent simple calculations on the finite field Z5 [13]. It is known that in that case, multiplication of all elements of the field with element of the field that is not zero, will lead to a permutation of the elements of the field [13]. Also, addition of the one nonzero element of the field to all other elements of the field will lead to a permutation of the elements of the field [13]. From that we can conclude that exactly one term will be congruent to 0 by modulo 5, and only one term will be congruent to 3 by modulo 5. That means that out of 5 threads defined by (B.1), three are going to stay after second step in Stage 2, which means that ¾ of the numbers that were left after step 1 in Stage 2, are going to stay after removal of the corresponding thread defined by number 5 (that is based on the Dirichelt's theorem [9] - all feasible threads defined by number 15 contain the same number of prime numbers).

After step 2 in Stage 2, all numbers can be written in the following forms

$$15n + 2,\ 15n + 11 \text{ and } 15n + 14. \qquad (B.3)$$

Now, this analysis can be applied on all consecutive step of Stage 2. In the step 3 of the Stage 2, we are going to apply a similar analysis like in the step 2 of Stage 2. In this case, instead of one thread defined by $3k + 2$, we have three threads defined by (B.3). In the third step of Stage 2, thread defined by number 7 is going to be removed. Impact of that removal is the easiest if we analyze the following forms of the remaining threads (here we are going to present forms for thread $15k + 2$; the other 2 threads could analyzed analogously)

$15(7n+1)+2$, $15(7n+2)+2$, $15(7n+3)+2$, $15(7n+4)+2$, $15(7n+5)+2$, $15(7n+6)+2$, $15(7n+7)+2$, or

$$105n+15*1+2,\ 105n+15*2+2,\ 105n+15*3+2,\ 105n+15*4+2,\ 105n+15*5+2,\ 105n+15*6+2,\ 105n+15*7+15*0+2. \quad \text{(B.4)}$$

Since all forms in (B.4) contain the term(s) divisible by 105 (and consequently divisible by 7), it is clear that additional forms that are going to be removed, will be removed based on the analysis of the following expressions

$$15*1+2,\ 15*2+2,\ 15*3+2,\ 15*4+2,\ 15*5+2,\ 15*6+2,\ 15*0+2, \quad \text{(B.5)}$$

or having in mind that $a*b \pmod 7 = a \bmod (7) * b \bmod (7)$, the forms of interest are given by

$$1*1+2,\ 1*2+2,\ 1*3+2,\ 1*4+2,\ 1*5+2,\ 1*6+2,\ 1*0+2. \quad \text{(B.6)}$$

Similarly to the situation in step 2, we can see that all seven terms in equation (B.6) represent simple calculations on the finite field Z7 [13]. Using the same line of reasoning like in the previous step, we can conclude that fraction of number of numbers that are going to stay after step 3 is exactly the one given in Table 2, and that is 5/6 of all numbers left after step 2 (here is assumed that the same analysis can be analogously performed for the other 2 threads defined by (B.3)). After this step 15 threads defined by number 105 are going to stay and each is going to contain the same percentage of prime numbers.

Now, it is obvious that proposed analysis can be applied to all consecutive steps of stage 2. In all cases, the removal of certain threads will be based on multiplication and addition of the finite field *Zpk*, where *pk* represents the odd prime number that defines thread that is going to be removed in the *k*-th step of Stage 2. In all cases those multiplications and addition will result in the permutation of all elements of the corresponding finite field and in every step, and it can be shown that in every step they are going to leave the ratio (*pk*-2)/(*pk*-1) of available numbers by using reasoning similar to the cases *pk* = {3, 5, 7}. From this analysis we can understand that the values in the second column of Table 2 are correct. Same can be concluded for all other threads that are not presented in the table. The proposed analysis holds also in the case of threads that are defined by prime numbers that are infinite (see [14]).

**APPENDIX C.**

Here, asymptotic density of numbers left, after implementation of the eSuS and the SGP sieve is going to be calculated (the eSuS represents extended $SuS^0$ sieve, in which after the removal of thread defined by some prime number in (2), also that prime number is removed). After removal of all even composite numbers and first *k* steps of the $SuS^0$ sieve, density of numbers left is given by the following equation

$$c_k = \frac{1}{2} \prod_{j=2}^{k+1} \left(1 - \frac{1}{p(j)}\right),$$

where $p(j)$ is *j*-th prime number. (We should have in mind that the density considered in previous equation are asymptotically correct).

In the case of $SuS^{-1}$ sieve the density of numbers left after the first *k*-steps is given by the following equation

$$c2_k = \prod_{j=2}^{k+1} \left(1 - \frac{1}{p(j) - 1}\right) = \prod_{j=2}^{k+1} \left(\frac{p(j) - 2}{p(j) - 1}\right).$$

So, if implementation of first sieve will result in the number of prime numbers smaller than *n* which we denote as $\pi(n)$, than implementation of the second sieve on some set of size $\pi(n)$ should result in the number of numbers $gp(n)$ that are defined by the following equation (for some big enough *n*)

$$gp(n) = r_{S2S1}(n) \cdot \frac{\pi(n)}{\ln(\pi(n))} \quad ,$$

where $r_{S2S1}(n)$ is defined by the following equation (*k* is the number of primes smaller or equal to $\sqrt{n}$)

$$r_{S2S1}(n) = \frac{c2_k}{c_k} = \frac{\prod_{2 < p \leq \sqrt{n}} \left(\frac{p - 2}{p - 1}\right)}{\prod_{p \leq \sqrt{n}} \left(\frac{p - 1}{p}\right)} = 2 \prod_{2 < p \leq \sqrt{n}} \left(\frac{p - 2}{p - 1}\right)\left(\frac{p}{p - 1}\right) \approx 2C_2 .$$

where *p* represents prime number. For *n* that is not big, $gp(n)$ should be defined as

$$gp(n) = f_{COR}(n) \cdot 2C_2 \cdot \frac{\pi(n)}{\ln(\pi(n))},$$

where $f_{COR}(n)$ represents correction factor that asymptotically tends toward 1 when $n$ tends to infinity. Here, the number of numbers left after eSuS is approximated by $n/\ln(n)$.